\documentclass[reqno]{amsart}

\input xy
\xyoption{all}
\usepackage{epsfig}
\usepackage{color}
\usepackage{amsthm}
\usepackage{amssymb}
\usepackage{amsmath}
\usepackage{amscd}
\usepackage{amsopn}
\usepackage{url}
\usepackage{hyperref}\hypersetup{colorlinks}

\usepackage{color} 

\definecolor{darkred}{rgb}{1,0,0} 
\definecolor{darkgreen}{rgb}{0,0.8,0}
\definecolor{darkblue}{rgb}{0,0,1}

\hypersetup{colorlinks,
linkcolor=darkblue,
filecolor=darkgreen,
urlcolor=darkred,
citecolor=darkgreen}

\newcommand{\labell}[1] {\label{#1}}

\numberwithin{equation}{section}
\newtheorem {Theorem}{Theorem}
\numberwithin{Theorem}{section}

\newtheorem {Proposition}[Theorem]{Proposition}
\newtheorem {Corollary}[Theorem]{Corollary}
\theoremstyle{definition}

\theoremstyle{remark}
\newtheorem{Remark}[Theorem]{Remark}
\newtheorem{Example}[Theorem]{Example}

\expandafter\chardef\csname pre amssym.def at\endcsname=\the\catcode`\@
\catcode`\@=11
\def\undefine#1{\let#1\undefined}
\def\newsymbol#1#2#3#4#5{\let\next@\relax
 \ifnum#2=\@ne\let\next@\msafam@\else
 \ifnum#2=\tw@\let\next@\msbfam@\fi\fi
 \mathchardef#1="#3\next@#4#5}
\def\mathhexbox@#1#2#3{\relax
 \ifmmode\mathpalette{}{\m@th\mathchar"#1#2#3}%
 \else\leavevmode\hbox{$\m@th\mathchar"#1#2#3$}\fi}
\def\hexnumber@#1{\ifcase#1 0\or 1\or 2\or 3\or 4\or 5\or 6\or 7\or 8\or
 9\or A\or B\or C\or D\or E\or F\fi}

\font\teneufm=eufm10
\font\seveneufm=eufm7
\font\fiveeufm=eufm5
\newfam\eufmfam
\textfont\eufmfam=\teneufm
\scriptfont\eufmfam=\seveneufm
\scriptscriptfont\eufmfam=\fiveeufm

\catcode`\@=\csname pre amssym.def at\endcsname

\newcommand{\CA}{{\mathcal A}}

\newcommand{\CN}{{\mathcal N}}

\newcommand{\Ham}{{\mathit Ham}}

\newcommand{\id}{{\mathit id}}

\newcommand{\tA}{\tilde{\mathcal A}}

\newcommand{\CB}{{\mathcal B}}
\newcommand{\Cc}{{\mathcal C}}

\def    \R      {{\mathbb R}}

\def    \Z      {{\mathbb Z}}

\def    \T      {{\mathbb T}}
\def    \CP     {{\mathbb C}{\mathbb P}}

\def    \12    {{\frac{1}{2}}}

\def    \p      {\partial}

\def    \HC     {\operatorname{HC}}

\def    \HQ     {\operatorname{HQ}}

\def    \H     {\operatorname{H}}

\def    \CF     {\operatorname{CF}}

\def    \Ham     {\operatorname{Ham}}

\def    \bx     {\bar{x}}

\def    \va     {\vec{a}}

\def    \MUCZ  {\operatorname{\mu_{\scriptscriptstyle{CZ}}}}

\def \odd   {\scriptscriptstyle{odd}}

\begin{document}


\setlength{\smallskipamount}{6pt}
\setlength{\medskipamount}{10pt}
\setlength{\bigskipamount}{16pt}





\title[Generic existence of periodic orbits]{On the generic existence of periodic orbits
  in Hamiltonian dynamics}

\author[Viktor Ginzburg]{Viktor L. Ginzburg}
\author[Ba\c sak G\"urel]{Ba\c sak Z. G\"urel}

\address{VG: Department of Mathematics, UC Santa Cruz,
Santa Cruz, CA 95064, USA}
\email{ginzburg@math.ucsc.edu}

\address{BG: Department of Mathematics, Vanderbilt University,
Nashville, TN 37240, USA} \email{basak.gurel@vanderbilt.edu}

\date{\today} \thanks{The work is partially supported by the NSF and
by the faculty research funds of the University of California, Santa
Cruz.}


\begin{abstract} 
  We prove several generic existence results for infinitely many
  periodic orbits of Hamiltonian diffeomorphisms or Reeb flows.  For
  instance, we show that a Hamiltonian diffeomorphism of a complex
  projective space or Grassmannian generically has infinitely many
  periodic orbits.  We also consider symplectomorphisms of the
  two-torus with irrational flux.  We show that such a
  symplectomorphism necessarily has infinitely many periodic orbits
  whenever it has one and all periodic points are non-degenerate.
\end{abstract}

\maketitle

\tableofcontents

\section{Introduction and main results}
\labell{sec:main-results}

\subsection{Introduction}
\label{sec:intro}
This paper focuses on the problem of $C^\infty$-generic existence of
infinitely many periodic orbits for Hamiltonian (or symplectic)
diffeomorphisms and Reeb flows. The topology of the underlying
symplectic or contact manifold plays an essential role in this
question and we prove several generic existence results which apply,
among other manifolds, to complex projective spaces and Grassmannians.
To put these results into perspective, recall that for a broad class
of closed symplectic manifolds, including all symplectically
aspherical ones, every Hamiltonian diffeomorphism has infinitely many
periodic orbits; see \cite{FrHa,Hi,Gi06conley,GG08}. However, in
contrast with the symplectically aspherical manifolds, complex
projective spaces and Grassmannians admit Hamiltonian diffeomorphisms
with only finitely many periodic orbits, and hence the existence of
infinitely many periodic orbits for these manifolds can be expected to
only hold $C^\infty$-generically or under suitable additional
requirements on the diffeomorphisms.

On the contact side, our results apply to many classes
of ``fillable'' contact forms. These include, among other examples,
various sets of contact forms on the unit cotangent bundles and
spheres, equipped with standard contact structures.  In particular, we
recover the $C^\infty$-generic existence of infinitely many closed
characteristics on convex hypersurfaces in $\R^{2n}$, \cite{Ek}, and
the $C^\infty$-generic existence of infinitely many closed geodesics,
\cite{Hi84,Ra89,Ra94}.

In both the Hamiltonian and contact cases, the proof is based on the
fact, established in \cite{GG08,GK09}, that indices and/or actions of
periodic orbits must satisfy certain relations when the
diffeomorphisms or the flow in question has only finitely many
periodic orbits. Roughly speaking, the argument is that these
relations are fragile and can be destroyed by a $C^\infty$-small
perturbation, and hence such a perturbation must create infinitely
many periodic orbits. (In the contact case, the resonance relations
that we utilize generalize those from \cite{Ek,EH,Vi89tr}.  This
reasoning is akin to the argument used in \cite{Ek,Ra94}.)

It is worth pointing out that a similar generic existence result
(Proposition \ref{prop:elliptic}) holds when $\H_{\odd}(M;\Z)\neq 0$
and is an easy consequence of the Birkhoff--Moser fixed point theorem,
\cite{Mo77}, and Floer's theory; see Section
\ref{sec:discussion}. However, in this case, there are no examples of
Hamiltonian diffeomorphisms with finitely many periodic orbits.

Finally, we also consider symplectomorphisms $\varphi\colon
\T^2\to\T^2$ with irrational flux or, equivalently, Hamiltonian
perturbations of an irrational shift of $\T^2$. We show that $\varphi$
necessarily has infinitely many periodic orbits whenever it has one
and all its periodic points are non-degenerate. The proof of this
theorem relies on the theory of Floer--Novikov homology developed in
\cite{LO,On}.

\begin{Remark}[Smoothness]
\label{rmk:smoothness}
Throughout the paper, for the sake of simplicity, all maps and vector
fields are assumed to be $C^\infty$-smooth and the spaces of maps or
vector fields are equipped with the $C^\infty$-topology, unless
specified otherwise. However, in all our results $C^\infty$ can be
replaced by $C^k$ with $k\geq 2$. The only exception is Proposition
\ref{prop:elliptic} where one has to require that $k\geq 4$.
\end{Remark}

\subsection{Periodic orbits of Hamiltonian diffeomorphisms}
\label{sec:ham}
Consider a closed symplectic manifold $(M^{2n},\omega)$, which
throughout this paper is assumed to be weakly monotone; see, e.g.,
\cite{HS} or \cite{MS04} for the definition. (This condition can be
eliminated by utilizing the machinery of virtual cycles.)  Recall that
$M$ is said to be \emph{monotone (negative monotone)} if
$[\omega]\mid_{\pi_2(M)}=\lambda c_1(M)\!\mid_{\pi_2(M)}$ for some
non-negative (respectively, negative) constant $\lambda$ and $M$ is called
\emph{rational} if $\left<[\omega], {\pi_2(M)}\right>=\lambda_0\Z$,
i.e., the integrals of $\omega$ over spheres in $M$ form a discrete
subgroup of $\R$. (When $\left<[\omega], {\pi_2(M)}\right>=0$, we set
$\lambda_0=\infty$.)  The constants $\lambda$ and $\lambda_0\geq 0$
are referred to as the \emph{monotonicity} and \emph{rationality
  constants}.  The positive generator $N$ of the discrete subgroup
$\left<c_1(M),\pi_2(M)\right>\subset \R$ is called the \emph{minimal
  Chern number} of $M$. When this subgroup is zero, we set $N=\infty$.
The manifold $M$ is called \emph{symplectically aspherical} if
$c_1(M)\mid_{\pi_2(M)} =0=[\omega]\mid_{\pi_2(M)}$. A symplectically
aspherical manifold is monotone and a monotone or negative monotone
manifold is rational.

A one-periodic (in time) Hamiltonian $H\colon \R/\Z\times M\to \R$
determines a time-dependent vector field $X_H$ on $M$ via Hamilton's
equation $i_{X_H} \omega = -dH$.  Let $\varphi=\varphi_H$ be the
Hamiltonian diffeomorphism of $M$, given as the time-one map of $X_H$.
Recall that there is a one-to-one correspondence between $k$-periodic
points of $\varphi$ and $k$-periodic orbits of $H$. In this paper, we
restrict our attention exclusively to periodic points of $\varphi$
such that the corresponding periodic orbits of $H$ are contractible. A
fixed point $x$ of $\varphi$ is said to be \emph{non-degenerate} if
$D\varphi_x\colon T_xM\to T_xM$ has no eigenvalues equal to one.
Recall also that $\varphi$ is called non-degenerate if all its fixed
points are non-degenerate. When all periodic points of $\varphi$ are
non-degenerate we will refer to $\varphi$ as a \emph{strongly
  non-degenerate} Hamiltonian diffeomorphism. Clearly, $\varphi$ is
strongly non-degenerate if and only if all iterations $\varphi^k$,
$k\geq 1$, are non-degenerate. In what follows, we will denote the
group of Hamiltonian diffeomorphisms of $M$, equipped with the
$C^\infty$-topology, by $\Ham(M,\omega)$.

Our first result concerns the generic existence of infinitely many
periodic orbits on symplectic manifolds with large $N$.

\begin{Theorem}
\labell{thm:ham1}
Assume that $n+1\leq N<\infty$. Then strongly non-degenerate
Hamiltonian diffeomorphisms with infinitely many periodic orbits form
a $C^\infty$-residual set in $\Ham(M,\omega)$.
\end{Theorem}

\begin{Example}
\labell{ex:cpn}
The only known monotone manifold to which this theorem applies is
$\CP^n$.   This manifold
admits Hamiltonian diffeomorphisms with finitely many fixed
points. The simplest of such diffeomorphisms is an irrational rotation
of $S^2$. Similar diffeomorphisms, arising from Hamiltonian torus actions,
exist in higher dimensions; see, e.g., \cite{GK09}.
These examples show that the genericity assumption in Theorem
\ref{thm:ham1} is essential. Note also that a Hamiltonian
diffeomorphism with finitely many periodic orbits need not be
associated with a Hamiltonian torus action. For instance, there exists
a Hamiltonian perturbation $\varphi$ of an irrational rotation of
$S^2$ with exactly three ergodic invariant measures: the Lebesgue
measure and the two measures corresponding to the fixed points of
$\varphi$; \cite{AK,FK}. 
There exist also a multitude of negative monotone manifolds satisfying
the hypotheses of the theorem. However, to the best of the authors'
knowledge there are neither examples of negative monotone manifolds admitting
Hamiltonian diffeomorphisms with finitely many periodic orbits nor any
results asserting that such manifolds do not exist.
\end{Example}

The theorem also holds when $N=\infty$, i.e.,
$c_1(TM)\mid_{\pi_2(M)}=0$, as can be easily seen by arguing as in
\cite{SZ}. However, in this case we expect the Conley conjecture to
hold, i.e., every Hamiltonian diffeomorphism to have infinitely many
periodic orbits. For instance, when, in addition, $M$ is rational,
this is proved in \cite{GG08}; see also
\cite{Gi06conley}. Furthermore, when $M=S^2$, a much stronger generic
existence result is established by the methods of two-dimensional dynamics in
\cite{We-H}. Namely, a $C^\infty$-generic Hamiltonian diffeomorphism
of $S^2$ has positive topological entropy and, as a consequence,
infinitely many hyperbolic periodic points.  (The proof relies on
Pixton's theorem asserting, roughly speaking, the generic existence of
hyperbolic periodic points of a particular type; see~\cite{Pi}.)

Theorem \ref{thm:ham1}, proved in Section \ref{sec:pf-ham}, is an easy
consequence of the mean index resonance relations established in
\cite{GK09}. A different type of relations, involving both the mean
indices and actions and proved in \cite{GG08}, leads to our next, more
technical, result. Denote by $\Lambda$ the Novikov ring of $M$,
equipped with the valuation $I_\omega(A):=-\left<\omega,A\right>$,
$A\in\pi_2(M)$, and by $*$ the pair-of-pants product in the quantum homology
$\HQ_*(M)$. (We refer the reader to, e.g., \cite{MS04} for a detailed
discussion of these notions; throughout this paper we adhere to
the conventions from \cite[Section 2]{GG08}.) Recall also that the Hofer
norm $\| \varphi\|$ of $\varphi\in \Ham(M,\omega)$ is defined as
$$
\|\varphi\|=\inf_H\int_0^1 (\max_M H_t-\min_M H_t)\, dt,
$$
where the infimum is taken over all $H$ such that $\varphi_H=\varphi$;
see, e.g., \cite{Po01}.

\begin{Theorem}
\labell{thm:ham2}
Assume that $M$ is monotone or negative monotone with monotonicity
constant $\lambda$ satisfying $|\lambda|<\infty$.
\begin{itemize}

\item[(i)] Then strongly non-degenerate Hamiltonian diffeomorphisms
with infinitely many periodic orbits form a $C^\infty$-residual set in
the Hofer ball 
$$
\CB=\{\varphi \in \Ham(M,\omega)\mid \|\varphi\|<\lambda_0\},
$$
where $\lambda_0$ is the rationality constant of $M$.

\item[(ii)] Assume in addition that there exists $u\in \H_{*<2n}(M)$ and $w\in
\H_{*<2n}(M)$ and $\alpha\in\Lambda$ such that 
\begin{equation}
\label{eq:hom}
[M]=(\alpha u)*w
\end{equation}
and one of the following requirements is satisfied 
\begin{itemize}

\item[(a)] $I_\omega(\alpha)=\lambda_0$;

\item[(b)] $2n-\deg u<2N$.
\end{itemize}
Then strongly non-degenerate Hamiltonian diffeomorphisms
with infinitely many periodic orbits form a $C^\infty$-residual set in
$\Ham(M,\omega)$.
\end{itemize}

\end{Theorem}

This theorem is proved in Section \ref{sec:pf-ham}.

\begin{Example} 
\label{ex:gap2}
Let us list some of the manifolds satisfying the hypotheses of
Theorem \ref{thm:ham2}(ii); see, e.g., \cite{GG08} and references therein
for details.
\begin{itemize}
\item The complex projective spaces $\CP^n$ and complex Grassmannians
  satisfy \eqref{eq:hom} and both (a) and (b).

\item Assume that $M$ satisfies \eqref{eq:hom} and (a) or (b) and $P$
  is symplectically aspherical. Then $M\times P$ satisfies
  \eqref{eq:hom} and (a) or, respectively, (b).

\item The product $M\times W$ of two rational manifolds satisfies
  \eqref{eq:hom} and (a) whenever $M$ does and
  $\lambda_0(W)=m\lambda_0(M)$, where $m$ is a positive integer or
  $\infty$. For instance, (a) holds for the products $\CP^n\times
  \CP^{m_1}\times\ldots \CP^{m_r}$ with $m_1+1,\ldots,m_r+1$ divisible
  by $n+1$ and equally normalized symplectic structures.

\item The monotone product $\CP^n\times W$, where $W$ is monotone and
  $\gcd\big(n+1,N(W)\big)\geq 2$, satisfies \eqref{eq:hom} and
  (b). For instance, this is the case for the monotone product
  $\CP^{n_1}\times\ldots\times\CP^{n_r}$ if
  $\gcd(n_1+1,\ldots,n_r+1)\geq 2$.

\end{itemize}
\end{Example}

While the proofs of Theorems \ref{thm:ham1} and \ref{thm:ham2} are
global and rely on relations between the mean indices and/or
actions of the periodic orbits, a local argument based on the
Birkhoff--Moser fixed point theorem (see \cite{Mo77}) combined with some
input from Floer's theory, yields the following.

\begin{Proposition}
\label{prop:elliptic}
Assume that $M$ is weakly monotone and $\H_{\odd}(M;\Z)\neq 0$. Then
strongly non-degenerate Hamiltonian diffeomorphisms with infinitely
many periodic orbits form a $C^\infty$-residual set in
$\Ham(M,\omega)$.
\end{Proposition}

This proposition, proved in Section \ref{sec:discussion}, also covers
the second case of Example \ref{ex:gap2}. Note that here the
requirement that $M$ be weakly monotone is purely technical and can
be eliminated completely with the use of virtual cycles.  It is also
worth pointing out that in the setting of the proposition, we have no
examples of Hamiltonian diffeomorphisms with finitely many periodic
points. It is possible that the Conley conjecture holds for symplectic
manifolds $M$ with $\H_{\odd}(M;\Z)\neq 0$, i.e., every Hamiltonian
diffeomorphism of $M$ has infinitely many periodic points.

Drawing on Theorems \ref{thm:ham1} and \ref{thm:ham2} and Proposition
\ref{prop:elliptic}, we conjecture that the existence of infinitely
many periodic points is a $C^\infty$-generic property, at least when
$M$ is monotone or negative monotone. Furthermore, it is perhaps
illuminating to look at these results in the context of the closing
lemma asserting, in particular, that the existence of a dense set of
periodic orbits is $C^1$-generic for both Hamiltonian diffeomorphisms
and flows; see \cite{PR}. Thus, once the $C^\infty$-topology is
replaced by the $C^1$-topology a much stronger result than the generic
existence of infinitely many periodic orbits holds -- the dense
existence. However, this is no longer true for the $C^k$-topology with
$k>\dim M$ as the results of M. Herman show (see \cite{He91a,He91b})
and the above conjecture on the $C^\infty$-generic existence of
infinitely many periodic orbits can be viewed as a viable form of a
$C^\infty$-closing lemma.

Furthermore, the aforementioned results of Herman (or rather the
argument from \cite{He91a,He91b}) also suggest that a $C^\infty$-small
Hamiltonian perturbation $\varphi$ of the shift $R_\theta$ of the standard
symplectic torus $\T^{2n}=\R^{2n}/\Z^{2n}$ never has periodic orbits
when $\theta$ satisfies a certain Diophantine condition.  (In
other words, here $R_\theta(x)=x+\theta$, where $x\in \T^{2n}$, and
$\varphi R^{-1}_\theta$ is a $C^\infty$-small Hamiltonian
diffeomorphism.) The situation however changes dramatically, as we
will see in the next section, when $\varphi$ is required to have at
least one periodic point, at least in the non-degenerate case and in
dimension two.

\subsection{Periodic orbits of symplectomorphisms of the two-torus} 
\label{sec:torus}
Let us consider a symplectomorphism $\varphi\colon \T^2\to\T^2$ such that
$\varphi R_\theta^{-1}$ is a Hamiltonian diffeomorphism for some
$\theta\in \T^2$. It is easy to see that this requirement is
equivalent to that the flux of $\varphi$ is $\theta$; see \cite{Ba}.

\begin{Theorem}
\labell{thm:torus}
Assume that $\varphi$ has at least one fixed point, all fixed points
of $\varphi$ are non-degenerate, and that at least one of the
components $\theta_1$ or $\theta_2$ of $\theta$ is irrational. Then
$\varphi$ has infinitely many periodic orbits. Moreover, for any
sufficiently large prime $k$, there is a simple $k$-periodic orbit.

\end{Theorem}

This theorem, proved in Section \ref{sec:pf-torus}, is a
symplectomorphism (two-dimensional) version of the non-degenerate
Conley conjecture established in \cite{SZ}; see also \cite{CZ86}.  As an
immediate consequence of Theorem \ref{thm:torus}, we obtain the
following result.

\begin{Corollary}
\label{cor:torus}
Assume that $\varphi$ is strongly non-degenerate and has at least one
periodic point, and that at least one of the components $\theta_1$ or
$\theta_2$ of $\theta$ is irrational. Then $\varphi$ has infinitely
many periodic orbits.
\end{Corollary}

\begin{Remark}
  The corollary still holds (without any non-degeneracy or fixed point
  requirement) when the flux $\theta$ is rational. Indeed, in this
  case a suitable iteration of $\varphi$ is Hamiltonian and one can
  apply the Conley conjecture, which in this case is proved in
  \cite{FrHa}.
\end{Remark}

One can view Corollary \ref{cor:torus} as a generic existence result
similar to Theorem \ref{thm:ham1} and \ref{thm:ham2}, with the class
of Hamiltonian diffeomorphisms replaced by symplectomorphisms with
flux $\theta$ and at least one homologically non-trivial periodic
point. However, this similarity is probably superficial and we
conjecture that any symplectomorphism of the standard $\T^{2n}$ with 
``sufficiently irrational'' flux and at least one homologically
essential periodic point has infinitely many periodic orbits. Note
that here, in contrast with many other Conley conjecture type results
(see, e.g., \cite{FrHa,Gi06conley,Hi} and references therein), the
assumption that the periodic point is homologically essential is
necessary as the following example due to John Franks indicates.

\begin{Example} Consider the one-form $\alpha_0=-\theta_2
  dq_1+\theta_1 dq_2$ on the two-torus $\T^2$ with angular coordinates
  $(q_1,q_2)$ and symplectic form $dq_1\wedge dq_2$.  The time-one map
  of the symplectic flow generated by $\alpha_0$ is $R_{\theta^0}$,
  where $\theta^0=(\theta_1,\theta_2)$. This map has no periodic
  points when at least one of the components of $\theta^0$ is
  irrational. Fix a small closed flow box $U$ in $\T^2$ and let
  $H$ be a smooth function on $U$ such that $H$ has only one critical
  point $x$ in $U$ and $dH=\alpha_0$ near $\p U$ and,
  finally, no Hamiltonian flow line of $H$ other than $x$ is entirely
  contained in $U$. Denote by $\alpha$ the closed
  one-form obtained from $\alpha_0$ by replacing $\alpha_0$ by $dH$ in
  $U$. It is easy to see that $x$ is the only periodic point of the
  time-one map $\varphi$ generated by $\alpha$. (By construction, the
  flow of $\alpha$ has no closed orbits other than the fixed point
  $x$.)  Furthermore, although the flux $\theta$ of $\varphi$ can be
  different from $\theta^0$, at least one component of $\theta$ is
  irrational by the Arnold conjecture. Thus, $\varphi$ satisfies all
  hypotheses of Theorem \ref{thm:torus} but the non-degeneracy. (As
  has been pointed out above, $\varphi$ is automatically a Hamiltonian
  perturbation $R_\theta$; see \cite{Ba}.) Note that the fixed point $x$ of
  $\varphi$ is degenerate and, moreover, $x$ can be destroyed by an
  arbitrarily small perturbation. In particular, the local Floer
  homology of $x$ is trivial; see \cite{Gi06conley,GG07gap}.
\end{Example}

\subsection{Periodic orbits of Reeb flows}
\label{sec:Reeb}
In this section we state the analogues of Theorems \ref{thm:ham1} and
\ref{thm:ham2} for Reeb flows. Similarly to the results from Section
\ref{sec:ham} that are based on the existence of resonance relations
for mean indices and actions of Hamiltonian diffeomorphisms, the
results of this section follow from the resonance relations for Reeb
flows, \cite{GK09}. 

To state these results, consider a closed contact manifold $(M,\xi)$.
Throughout this section we assume that $c_1(\xi)=0$ although this
condition may in some instances be relaxed. We refer the reader to
\cite{Ge} for a general introduction to contact topology; see also
\cite{Bo,BO,El}.  Having the contact manifold $(M,\xi)$ fixed, we only consider
contact forms $\alpha$ on $M$ with $\ker\alpha=\xi$.  We say that the
\emph{Weinstein conjecture holds for $(M,\xi)$ if the Reeb flow for
  every contact form $\alpha$ has a periodic orbit}; see, e.g.,
\cite{Ge}.  We call a contact form \emph{non-degenerate} when
all periodic orbits of its Reeb flow are non-degenerate. (This notion
is similar to strong non-degeneracy for Hamiltonian diffeomorphisms.)
For a non-degenerate closed Reeb orbit $x$, we set its degree to be
$|x|=\MUCZ(x)+n-3$, where $\MUCZ(x)$ stands for the Conley--Zehnder
index of $x$ (with its standard normalization), and denote by
$\Delta(x)$ the mean index of $x$.

Our generic existence results concern sets of contact structures
satisfying some additional conditions. Namely, let $\Cc$ be a set of contact
forms $\alpha$ (such that $\ker\alpha=\xi$), equipped with the
$C^\infty$-topology, meeting the following two requirements:

\begin{itemize}

\item[(C1)] The intersection of the set of non-degenerate contact forms with
$\Cc$ is a residual set in $\Cc$.

\item[(C2)] For every non-degenerate form $\alpha\in\Cc$, a periodic
  orbit $x$ of the Reeb flow of $\alpha$ and a small neighborhood $U$
  of $x$, there exists a sequence of perturbations $\alpha_k\to\alpha$
  in $\Cc$ such that $\alpha_k-\alpha$ is supported in $U$ and
  $\Delta(x_k)\neq \Delta(x)$ for all $k$, where $x_k$ is the periodic
  orbit of the Reeb flow of $\alpha_k$ arising as a perturbation
  of~$x$.

\end{itemize}
Note that the orbit $x_k$ in (C2) is unique and close to $x$ once $U$
and $\alpha-\alpha_k$ are sufficiently small and that
$\Delta(x_k)\to\Delta(x)$ due to continuity of $\Delta$; see
\cite{SZ}. Thus, condition (C2) is the requirement that
$\Delta(x)$ can be varied by varying $\alpha$ in $\Cc$.

To state the first result, assume that $M$ is the boundary of a
compact manifold $W$. We say that $(M,\xi)$ is \emph{symplectically
  fillable} by $W$ if $d\alpha_0$, for some contact form $\alpha_0$ on
$(M,\xi)$, extends to a symplectically aspherical symplectic structure
$\omega_0$ on $W$. Then, as is not hard to see, the same is true for
every contact form $\alpha$ on $(M,\xi)$. The extension $\omega$ of $d\alpha$ is
obtained from $\omega_0$ by modifying the latter near $M=\p W$, and
$(M,\alpha_0)$ and $(M,\alpha)$ have the same linearized contact
homology with respect to these fillings.

\begin{Theorem}
\label{thm:cont2}
Assume that $(M,\xi)$ is symplectically fillable, the Weinstein
conjecture holds for $\xi$, and that $\Cc$ satisfies (C1) and (C2).
Then non-degenerate contact forms with infinitely many periodic orbits
form a residual subset in $\Cc$.
\end{Theorem}

\begin{Example}
  Assume that $(M,\xi)$ is symplectically fillable. Then requirements
  (C1) and (C2) are obviously satisfied for the set of all contact
  forms $\alpha$ on $(M,\xi)$.
\end{Example}

The second result concerns the situation where the contact manifold $(M,\xi)$ 
is not required to be fillable.

\begin{Theorem}
\label{thm:cont1}
Assume that the Weinstein conjecture holds for $\xi$, the set
$\Cc$ meets requirements (C1), (C2) and also the following:
\begin{itemize}
\item[(C3)] the Reeb flow of any non-degenerate contact form in $\Cc$ 
 has no contractible periodic orbits $x$ with $|x|=0$ or $\pm 1$. 
\end{itemize}
Then non-degenerate contact forms with infinitely many periodic orbits
form a residual subset in $\Cc$.
\end{Theorem}

\begin{Remark}
  Requirement (C3) is quite restrictive: there exist
  contact structures which admit no non-degenerate contact forms
  satisfying (C3). For instance, this is the case
  for an overtwisted contact structure on $S^3$; see, e.g., 
  \cite{Yau}
\end{Remark}

\begin{Example}
\label{ex:geodesics}
Let $M$ be the unit cotangent bundle $ST^*P$ of a closed manifold $P$,
equipped with the standard contact structure and let $\Cc$ be the set
of contact forms on $M$ associated with Riemannian metrics on
$P$. Then $(M,\xi)$ is fillable and $\Cc$ satisfies requirements (C1)
and (C2); condition (C3) is met whenever $\dim P >3$. (See
\cite{Ab,An} for the proof of (C1) and \cite{AS04,Lo,Vi90} and
references therein and, in particular \cite{Du}, for a discussion of
indices; (C2) can be established similarly to the argument from
\cite{KT}.)  Likewise, let $M=S^{2n-1}$ be equipped with the standard
contact structure and let $\Cc$ be the class of contact forms arising
from embeddings of $M$ into $\R^{2n}$ as a strictly convex
hypersurface enclosing the origin. Then (C1) is satisfied and it is
not hard to see that $\Cc$ also meets requirements (C2) and (C3); see,
e.g., \cite{Lo} and references therein.  Furthermore, any contact form
$\alpha$ on $S^{2n-1}$ giving rise to the standard contact structure
is symplectically fillable by $W=B^{2n}$.  (Indeed, $\alpha$ is the
restriction of the form $\sum (x_i dy_i-y_i dx_i)/2$ on $\R^{2n}$ to
an embedding $S^{2n-1}\hookrightarrow \R^{2n}$ bounding a starshaped
domain.)  Thus, Theorems \ref{thm:cont2} and \ref{thm:cont1}
generalize the $C^\infty$-generic existence of infinitely many closed
characteristics on convex hypersurfaces in $\R^{2n}$ (see \cite{Ek})
and the $C^\infty$-generic existence of infinitely many closed
geodesics (see \cite{Hi84,Ra89,Ra94}). Similar considerations apply,
of course, to Finsler metrics, symmetric as well as asymmetric.

\end{Example}

\begin{Remark} In Theorems \ref{thm:cont2} and \ref{thm:cont1} one can
  also require the homotopy classes of closed orbits to lie in a fixed
  set of free homotopy classes of loops in $M$, closed under
  iterations. (For instance, one can require the orbits to be
  contractible.)  However, in this case the Weinstein conjecture must
  also hold for such orbits, which is a non-trivial restriction on the
  contact structure and the set of homotopy classes. (The argument
  from Section \ref{sec:pf-reeb} readily proves this generalization of
  the theorems, cf.\ the discussion in \cite[Section 1.3]{GK09}.)
\end{Remark}

\subsection{Acknowledgments} We are grateful to Alberto
Abbondandolo, Christian Bonatti, Yasha Eliashberg, John Franks,
Anatole Katok, Leonid Polterovich, and Marcelo Viana for useful
comments and remarks.

\section{Proofs and remarks}

\subsection{Hamiltonian diffeomorphisms: Proofs of Theorems
  \ref{thm:ham1} and \ref{thm:ham2}}
\labell{sec:pf-ham}

\begin{proof}[Proof of Theorem \ref{thm:ham1}]
  Recall that to a contractible periodic orbit $x$ of $H$ we can
  associate the mean index $\Delta(x)\in \R/2N\Z$ as in \cite{SZ}.
  Strictly speaking, the mean index, viewed as a real number, depends
  on the choice of the capping of $x$. However, $\Delta(x)$ is well
  defined once it is regarded as an element of the circle $\R/2N\Z$.
  Consider the set $\Delta^\infty$ formed by the indices of simple
  contractible periodic orbits of $H$ for all periods. (Here we treat
  $\Delta^\infty$ as a genuine set: if two orbits
  have equal indices, their index enters the collection only once.)
  Denote by $\Delta^{k}$ the subset of $\Delta^\infty$ formed by the
  mean indices of periodic orbits with period less than 
  $k$. Note that the set $\Delta^{k}$ is necessarily finite when
  $H$ is strongly non-degenerate.

  Recall also that, as readily follows from 
\cite[Theorem 1.1 and Remark 1.6]{GK09},
  whenever $\varphi_H$ has finitely many periodic orbits, the set
  $\Delta^\infty=(\Delta_1,\ldots,\Delta_m)$ satisfies a resonance relation
  of the form
\begin{equation}
\label{eq:ham-res}
a_1\Delta_1+\ldots+a_m\Delta_m=0\mod 2N
\end{equation}
for some non-zero vector $\va=(a_1,\ldots,a_m)\in\Z^m$.

Let $\CN_k(\va)$ be the set of strongly non-degenerate Hamiltonian
diffeomorphisms $\varphi$ that do not satisfy the resonance relation
\eqref{eq:ham-res} up to period $k$. (Here we do not require $\varphi$ to
have finitely many periodic points.) More precisely, a
non-degenerate Hamiltonian diffeomorphism is in $\CN_k(\va)$ if either
the number of non-zero components in the vector $\va$ exceeds the
cardinality of $\Delta^{k}$ or \eqref{eq:ham-res} fails for any
choice of a subset $\Delta=(\Delta_1,\ldots,\Delta_m)$ in $\Delta^{k}$
and any ordering of this subset. 

It is routine to show that $\CN_k(\va)$ is an open, dense subset in
$\Ham(M,\omega)$ in the $C^\infty$-topology. Indeed, it is clear that
$\CN_k(\va)$ is open. To prove that it is dense, let us consider a
non-degenerate Hamiltonian diffeomorphism satisfying the
resonance relation \eqref{eq:ham-res} for some subset $\Delta\subset
\Delta^{k}$ formed by $\Delta_i=\Delta(x_i)$ where all orbits
$x_i$ have period less than or equal to $k$ and, say, $a_1\neq 0$. (We
will assume for the sake of simplicity that $\Delta$ is the only
ordered collection in $\Delta^{k}$ satisfying
\eqref{eq:ham-res} and that $x_1$ is the only orbit with mean index
$\Delta_1$ and period less than $k$. The
general case can be dealt with in a similar fashion.)  By applying a sufficiently
$C^\infty$-small perturbation of $H$ localized near $x_1$ in the
``space-time'' $S^1\times M$, we can change the value of $\Delta(x_1)$
without affecting other periodic orbits of period up to $k$. This
change will destroy the resonance relation \eqref{eq:ham-res}. It will
also create no periodic orbits of period less than $k$
(since $x_1$ is non-degenerate), and hence no new resonance relations
for such orbits. Thus, $\CN_k(\va)$ is dense.

Taking the intersection of these sets for all
$k$ and $\va$, we obtain a $C^\infty$-residual
subset $\CN$ of $\Ham(M,\omega)$, which, by the result from \cite{GK09}
quoted above, contains no Hamiltonian diffeomorphisms with finitely
many periodic orbits.
\end{proof}

\begin{proof}[Proof of Theorem \ref{thm:ham2}] The argument is similar
  to the proof of Theorem \ref{thm:ham1}, but instead of the resonance
  relations \eqref{eq:ham-res} from \cite{GK09} we utilize a relation
  involving both the mean indices and actions and proved in
  \cite{GG08}.

  Let, as above, $H$ be a one-periodic in time Hamiltonian on $M$ and
  let $x$ be a $k$-periodic orbit of $H$.  The normalized augmented action of $H$
  on $x$ is defined as
$$
\tA_H(x)=\big(\CA_H(\bx)-\lambda \Delta_H(\bx)\big)/k, 
$$
where $\bx$ is the orbit $x$ equipped with an arbitrary capping.  Here
$\CA_H(\bx)$ stands for the ordinary action of $H$ on $\bx$ and
$\Delta_H(\bx)$ is mean index of $H$ on $\bx$.  It is clear that
$\tA_H(x)$ is independent of the capping.  This definition, borrowed
from \cite{GG08}, is inspired by the considerations in \cite[Section
1.6]{Sa} and \cite[Section 1.4]{EP}, where the Conley--Zehnder index
is utilized in place of the mean index. For us, the main advantage of
using the mean index is that iterating an orbit does not change the
normalized augmented action, i.e., $\tA_H(x^l)=\tA_H(x)$.  The reason is
that the ordinary action and the mean index are both homogeneous:
$\CA_H(\bx^l)=l\CA_H(\bx)$ and $\Delta_H(\bx^l)=l\Delta_H(\bx)$.
Moreover, geometrically identical orbits have equal normalized
augmented action. (Two periodic orbits of $H$ are said to be
geometrically identical if the corresponding periodic orbits of
$\varphi_H$ coincide as subsets of $M$; see \cite[Section~1.3]{GG08}.)

By \cite[Corollary 1.11]{GG08}, under the hypotheses of the theorem,
there exist two geometrically distinct periodic orbits $x$ and $y$ of
$H$ with $\tA_H(x)=\tA_H(y)$ whenever $\varphi_H$ has finitely many
periodic orbits. In case (ii), let $\CN_k$ be the set of strongly
non-degenerate Hamiltonian diffeomorphisms $\varphi_H$ such that all
geometrically distinct periodic orbits of $\varphi_H$ up to period $k$
have different normalized augmented actions.  As in the proof of
Theorem \ref{thm:ham1}, it is easy to show that $\CN_k$ is open and
dense in $\Ham(M,\omega)$ in the $C^\infty$-topology. (In case (i),
$\CN_k$ is defined similarly, but only as a subset of $\CB$. Then
$\CN_k$ is open and dense in $\CB$.)  The set $\CN=\bigcap_{k}\CN_k$
is residual and contains no Hamiltonian diffeomorphisms with finitely
many periodic orbits.
\end{proof}

\subsection{Symplectomorphisms of $\T^2$: Proof of Theorem \ref{thm:torus}}
\labell{sec:pf-torus}
The proof of the theorem relies on the machinery of Floer--Novikov
homology developed in \cite{LO,On} and throughout the proof we 
use the conventions and notation from these works.

\begin{proof}[Proof of Theorem \ref{thm:torus}] 
  Let us assume first that $\varphi^k$ is a non-degenerate iteration
  of $\varphi$ and denote by $\CF_*(\varphi^k)$ the Floer--Novikov
  complex of $\varphi^k$. (Strictly speaking, to define this complex,
  we need to pick a path $\psi_t$ connecting $\varphi^k$ to the
  identity. To this end, we fix a path $\varphi_t$ from
  $\varphi$ to $\id$ and set $\psi_t=\varphi_t^k$. The choice of the
  path $\varphi_t$ is immaterial.) The complex $\CF_*(\varphi^k)$ is
  generated by the fixed points of $\varphi^k$, corresponding to the
  contractible $k$-periodic orbits of $\varphi_t$, over the Novikov
  ring $\Lambda_{k\theta}$; see \cite{LO}. It is essential for what follows that the
  elements of $\Lambda_{k\theta}$ have zero degree. (To be quite
  precise, one has to consider the lifts of the fixed points to the
  universal covering $\R^2$ of $\T^2$. Note also that the Novikov ring
  $\Lambda_{k\theta}$ depends on $k$ if one of the components of
  $\theta$ is rational.)  

The complex $\CF_*(\varphi^k)$ is acyclic.  Indeed, $\varphi^k$ is a
Hamiltonian deformation of $R_{k\theta}$, and hence the two
symplectomorphisms have the same Floer--Novikov homology. Furthermore,
since at least one of the components of $\theta$ is irrational
$R_{k\theta}$ has no fixed points, and hence its Floer--Novikov homology vanishes.

Arguing by contradiction, assume now that $\varphi$ has no simple
$k$-periodic points for a (large) prime $k$. Then the fixed points of
$\varphi^k$ are the $k$th iterations of the fixed points of
$\varphi$. (Here and in what follows, we consider only the fixed points
corresponding to the contractible orbits of $\varphi_t$.) Furthermore,
once $k$ is greater than the degree of any root of unity among the
Floquet multipliers of the one-periodic orbits of $\varphi_t$, the
iteration $\varphi^k$ is non-degenerate and the complex
$\CF_*(\varphi^k)$ is generated by the $k$th iterations of the fixed
points of $\varphi$.  Let us group these fixed points according to
their mean indices, placing all fixed points with the same index into
one group. Thus, we have $r$ groups corresponding to different real
numbers $\Delta_1,\ldots,\Delta_r$ occurring as the mean indices for
$\varphi$. This grouping is carried over to the $k$th iteration
$\varphi^k$. When $k$ is so large that $k|\Delta_i-\Delta_j|>3$ (if
$i\neq j$), the complex $\CF_*(\varphi^k)$ breaks down into a direct
sum of complexes each of which is generated by the $k$th iteration of
the fixed points from one group.  Each of these complexes is acyclic.
Our goal is to show that this is impossible.

To this end, consider the complex generated by the orbits in one group,
say, $x_1^k,\ldots,x_m^k$. (Since $\varphi$ has a fixed point, there is
at least one non-empty group.) Since all orbits $x_i^k$ are
non-degenerate, we have
$$
|\MUCZ(x_i^k)-k\Delta|<1,
$$
where $\Delta$ is the mean index of the group; see \cite{CZ86,SZ}. 

Then $\MUCZ(x_i^k)=k\Delta$ if $k\Delta\in \Z$. Thus, in this case all
generators $x_i^k$ have the same degree, which is impossible, for the
complex is acyclic. When $k\Delta\not\in \Z$, every orbit $x_i^k$ is
necessarily elliptic: its Floquet multipliers $\lambda,\bar{\lambda}$
are on the unit circle and different from
$\pm 1$. For such an orbit, in dimension two, the mean index
completely determines the Conley--Zehnder index. Hence again, we
arrive at a contradiction with acyclicity, for all generators $x_i^k$
have the same degree.
\end{proof}

\subsection{Reeb flows: Proofs of Theorems \ref{thm:cont2} and \ref{thm:cont1}}
\label{sec:pf-reeb}
The proof of both of these results is quite similar to the proofs of
Theorems \ref{thm:ham1} and \ref{thm:ham2} except that now a different
resonance relation is used.

\begin{proof}[Proof of Theorems \ref{thm:cont2} and \ref{thm:cont1}]
  Let $\Cc$ be as in either of the theorems and let $\alpha\in\Cc$ be
  non-degenerate. Following \cite{GK09}, we call a simple periodic
  orbit $x$ of the Reeb flow of $\alpha$ \emph{bad} if the linearized
  Poincar\'e return map along $x$ has an odd number of real
  eigenvalues strictly smaller than $-1$. Otherwise, the orbit is said
  to be \emph{good}. (This terminology differs slightly from the
  standard usage, cf.\ \cite{Bo,BO}.)  When the orbit $x$ is good, the
  parity of the Conley--Zehnder indices $\MUCZ(x^k)$ is independent of
  $k$; if $x$ is bad, the parity of $\MUCZ(x^k)$ depends on the
  parity of $k$. We denote the mean index of an orbit $x$ by
  $\Delta(x)$ and set
  $\sigma(x)=(-1)^{|x|}=(-1)^{n+1}(-1)^{\MUCZ(x)}$. In other words,
  $\sigma(x)$ is, up to the factor $(-1)^{n+1}$, the topological index of
  the orbit $x$ or, more precisely, of the Poincar\'e return map along
  $x$.

  Furthermore, under the hypotheses of either of the theorems, the
  contact homology $\HC_*(M,\xi)$ of $(M,\xi)$ is defined and
  independent of $\alpha\in \Cc$; see \cite{Bo,BO}. (In the case of
  Theorem \ref{thm:cont2}, we use the fact, mentioned in Section
  \ref{sec:Reeb}, that the linearized contact homology is independent of $\alpha$
  as long as the fillings are adjusted accordingly.)  As
  in \cite{GK09}, set
 \begin{equation}
\label{eq:Euler1}
\chi^\pm(W,\xi)=\lim_{N\to\infty}\frac{1}{N}
\sum_{l=l_\pm}^N(-1)^l\dim \HC_{\pm l}(W,\xi),
\end{equation}
where $l_-=-2$ and $l_+=2n-4$, provided that all terms are finite and
the limits exist, and let
$$
\chi(W,\xi):=\frac{\chi^+(W,\xi)+\chi^-(W,\xi)}{2}.
$$
We call
$\chi(W,\xi)$ the mean Euler characteristic of $\xi$. (This invariant
is also considered in \cite[Section 11.1.3]{VK:thesis}.)

Then, whenever the Reeb flow of $\alpha$ has finitely many simple
periodic orbits, $\dim \HC_{\pm l}(W,\xi)<\infty$ when $\pm l>
l_{\pm}$, the limits in \eqref{eq:Euler1} exist, and the mean indices
of the orbits satisfy the resonance relation
\begin{equation}
\label{eq:Euler2}
{\sum} \frac{\sigma(x_i)}{\Delta(x_i)}
+\frac{1}{2}{\sum} \frac{\sigma(y_i)}{\Delta(y_i)}=
\chi(W,\xi),
\end{equation}
where the first sum is over all good simple periodic orbits $x_i$, the
second sum is over all bad simple periodic orbits $y_i$ and in both
cases the orbits with zero mean index are excluded.  (See \cite{GK09}
for a proof; this result generalizes the resonance relations from
\cite{EH,Vi89tr}.)

Denote by $\CN_k$ the set of non-degenerate forms $\alpha\in \Cc$
such that \eqref{eq:Euler2} fails when the summation on the left hand
side is taken over all simple periodic orbits of period less than
$k$. Note that, even though now we are not assuming that the Reeb flow
has finitely many periodic orbits, the number of orbits with 
period bounded from above by $k$ is finite as a consequence of non-degeneracy. 

We claim that, when $k$ is sufficiently large, the set $\CN_k$ is
$C^\infty$-open and dense in $\Cc$. Indeed, it is clear that $\CN_k$
is open.  To show that it is dense, it suffices, by (C1), to prove
that any neighborhood of a non-degenerate form $\alpha_0\in\Cc$
contains a form $\alpha\in\CN_k$. This is clear due to (C2): by
varying $\alpha_0$ in $\Cc$, one can make \eqref{eq:Euler2} fail for
orbits with period bounded by $k$. (At this point we use the
assumption that the Weinstein conjecture holds for $\xi$ to make sure
that the left hand side of \eqref{eq:Euler2} contains at least one
term when $k$ is large.)

Finally, taking the intersection of the sets $\CN_k$, for all large
$k$, we obtain a $C^\infty$-residual set $\CN\subset\Cc$ of
non-degenerate forms and, by \eqref{eq:Euler2}, the Reeb flow of any
form in $\CN$ has infinitely many periodic orbits.
\end{proof}

\subsection{Discussion} 
\label{sec:discussion}
One common requirement in the majority of generic existence results for infinitely
many periodic orbits is the existence of one such
orbit, which in some instances must satisfy certain additional conditions.

To illustrate this point, recall that the Birkhoff--Moser fixed point theorem (see
\cite{Mo77}) asserts that when $x$ is a non-degenerate, non-hyperbolic
periodic point of a Hamiltonian diffeomorphism $\varphi$, a
$C^\infty$-generic perturbation of $\varphi$ (which can be taken to be
strongly non-degenerate) has infinitely many periodic points in a
neighborhood of $x$. (A $k$-periodic point $x$ is called
non-hyperbolic if at least one of the eigenvalues of the linearized
map $d\varphi^k\colon T_xM\to T_xM$ is on the unit circle.) Thus, the
Birkhoff--Moser fixed point theorem implies the $C^\infty$-generic
existence of infinitely many periodic orbits in a neighborhood (in the
$C^\infty$-topology) of a Hamiltonian diffeomorphism with a
non-degenerate, non-hyperbolic periodic point.  Then one may try to
deal with the existence problem in the case where all periodic points
are hyperbolic by some different, usually \emph{ad hoc}, method. For
instance, this is the approach used in \cite{Hi84} to establish the
generic existence of infinitely many geodesics on simply connected,
compact symmetric spaces of rank one (e.g., spheres). In the case of Hamiltonian
diffeomorphisms, arguing along these lines we obtain Proposition
\ref{prop:elliptic}.

\begin{proof}[Proof of Proposition \ref{prop:elliptic}] By the
  Birkhoff--Moser fixed point theorem it suffices to show that a
  strongly non-degenerate Hamiltonian diffeomorphism $\varphi$ has a
  non-hyperbolic periodic point or $\varphi$ has infinitely many
  periodic orbits. Assume that all periodic points of $\varphi$ are
  hyperbolic. Then all periodic points of $\varphi$ have even
  Conley--Zehnder indices, which is impossible by Floer's theory since
  $\H_{\odd}(M;\Z)\neq 0$ (see, e.g., \cite{Sa,SZ}), or every
  iteration of $\varphi$ has a hyperbolic point with negative real
  eigenvalues. It follows then that every iteration of the form
  $\varphi^{2^k}$ must have a simple periodic orbit and thus the
  number of periodic orbits is infinite.
\end{proof}

In a similar vein, the existence of a periodic point (hyperbolic)
plays an important role in the argument from \cite{We-H}, mentioned in
Section \ref{sec:ham}, concerning Hamiltonian diffeomorphisms of
$S^2$.

The method employed in this paper and utilizing the resonance relations  is
no exception in that it also relies, implicitly or explicitly, on the
existence of one periodic point, even though we use global
symplectic-topological rather than dynamical systems arguments. In
Theorems \ref{thm:ham1} and \ref{thm:ham2}, the existence of one fixed
point is implicit, for this requirement is automatically satisfied due
to the Arnold conjecture. In Theorem \ref{thm:torus}, the requirement
is explicit and essential; see Section \ref{sec:torus}. In Theorems
\ref{thm:cont2} and \ref{thm:cont1}, the requirement is also explicit
(the Weinstein conjecture), but is conjecturally always met.

\begin{Remark} 
  We conclude this discussion by pointing out an aspect of the generic
  existence problem for periodic orbits that is not touched upon in
  this paper (except in Theorem \ref{thm:torus}). This is the question
  of the (generic) growth of the number of periodic orbits, which is
  of interest even when infinitely many periodic orbits exist
  unconditionally. For instance, when the manifold is symplectically
  aspherical, a generic Hamiltonian diffeomorphism has a simple
  periodic orbit for every sufficiently large prime period; see
  \cite{SZ}. Thus, the number of geometrically distinct periodic
  orbits of period less than or equal to $k$ generically grows at
  least as $k/\log k$. A similar lower bound
  holds, up to a factor, for closed geodesics (cf.\ Example
  \ref{ex:geodesics}) on a Riemannian manifold with a finite, but
  non-trivial, fundamental group, \cite{BTZ}. However, to the best of
  the authors' knowledge no such results for generic Hamiltonian
  diffeomorphisms of, say, $\CP^n$ or under the hypotheses of
  Proposition \ref{prop:elliptic} have been obtained.
\end{Remark}

\end{document}